\documentclass[12pt]{amsart}
\usepackage{amsmath,amsfonts,amssymb}
\usepackage{aliascnt}
\usepackage{aliascnt}
\usepackage{graphicx}
\usepackage{mathrsfs}
\usepackage{amscd}
\usepackage{marginnote}
\usepackage{color}

\usepackage[bookmarks=true,pdfstartview=FitH, pdfborder={0 0 0}, colorlinks=true,citecolor=blue, linkcolor=blue]{hyperref}

\newtheorem{theorem}{Theorem} %[section]
\newaliascnt{conj}{theorem}
\newaliascnt{cor}{theorem}
\newaliascnt{lemma}{theorem}
\newaliascnt{fact}{theorem}
\newaliascnt{claim}{theorem}
\newaliascnt{prop}{theorem}
\newaliascnt{definition}{theorem}
\newaliascnt{qn}{theorem}
\newaliascnt{assump}{theorem}

\newtheorem{cor}[cor]{Corollary}

\aliascntresetthe{conj} \aliascntresetthe{cor}
\aliascntresetthe{lemma} \aliascntresetthe{fact}
\aliascntresetthe{claim} \aliascntresetthe{prop}
\aliascntresetthe{definition} \aliascntresetthe{qn}
\aliascntresetthe{assump}

\theoremstyle{definition}
\newaliascnt{example}{theorem}

\aliascntresetthe{example}

\theoremstyle{remark}
\newaliascnt{rmk}{theorem}
\newtheorem{remark}[rmk]{Remark}
\aliascntresetthe{rmk} 

\def\sek~{\S{}}

\newcommand{\p}{\partial}
\newcommand{\F}{\mathcal{F}}

\makeatletter
\def\blfootnote{\xdef\@thefnmark{}\@footnotetext}
\makeatother

\begin{document}
\title{Non-existence of certain singularities in Legendrian foliations}

\author{Yang Huang}
\address{Department of Mathematics -- Centre for Quantum Geometryof Moduli Spaces, Ny Munkegade 118, Building 1530, 8000 Aarhus C, Denmark}
\email{yhuang@qgm.au.dk}

\blfootnote{Y.H was supported by the Center of Excellence Grant ``Centre for Quantum Geometry of Moduli Spaces'' from the Danish National Research Foundation (DNRF95).}

\begin{abstract}
In this paper we show that the singular locus of a Legendrian foliation as defined in \cite{H2013} is a compact submanifold whose connected components are of codimension at most two. As a consequence, given any closed $(n+1)$-dimensional coisotropic submanifold $W$ in a contact $(2n+1)$-manifold, the contact structure in a sufficiently small neighborhood of $W$ is uniquely determined by the characteristic (Legendrian) foliation.
\end{abstract}

\maketitle

Let $(M^{2n+1},\xi)$ be a contact manifold with contact structure $\xi=\ker\alpha$. Consider a closed $(n+1)$-submanifold $W \subset M$ together with the restricted one-form $\lambda=\alpha|_W$. We say $W$ is a {\em coisotropic submanifold} if $\lambda \wedge d\lambda=0$ on $W$. In other words, one can integrate $\ker\lambda$ into a singular foliation $\F$, which we will call the {\em characteristic foliation}. It is easy to check that the leaves of $\F$ are Legendrian submanifolds, so we sometimes also call $\F$ a {\em Legendrian foliation}.

We are interested in the singular locus of $\F$, which is defined to be $S=\{\lambda=0\}$. In a previous work \cite{H2013}, we see that each path-connected component of $S$ must belong to one of the following three kinds.
\begin{itemize}
	\item[$(A_0)$:] A closed orientable $(n-1)$-dimensional submanifold;
	\item[$(A_1)$:] An open orientable $(n-1)$-dimensional submanifold;
	\item[$(A_2)$:] A $C^1$-smooth compact $n$-dimensional submanifold.
\end{itemize}
Moreover we see that if $S$ consists of only singularities of type $(A_0)$ and $(A_2)$, then the germ of the contact structure in a neighborhood of $W$ is uniquely determined by $\F$. Conversely, given any singular foliation $\F$ with only singularities of type $(A_0)$ and $(A_2)$, one can construct a standard contact neighborhood of $W$ such that $W$ becomes a coisotropic submanifold with characteristic foliation equal to $\F$.

Using the above terminologies, we can formulate the main result of this paper as follows.

\begin{theorem} \label{thm:main}
Given a closed coisotropic submanifold $W^{n+1} \subset (M^{2n+1},\xi)$ with Legendrian foliation $\F$, there exists no singularities of type $(A_1)$ in $\F$.
\end{theorem}

Combined with the existence and uniqueness results in \cite{H2013}, we obtain the following result.

\begin{cor} \label{cor:main}
The contact germ in a neighborhood of a coisotropic submanifold is uniquely determined by the characteristic foliation.
\end{cor}

\begin{remark}
\autoref{cor:main} is well-known in three-dimensional contact geometry but the proof in that case involves only Moser's technique.
\end{remark}

\section{Singularities in Legendrian foliations} \label{sec:sing}

In this section, we recall some relevant general properties of singularities in Legendrian foliations from \cite{H2013}.

Let $(W,\F) \subset (M,\xi)$ be a coisotropic submanifold with Legendrian foliation $\F$. Suppose, for the rest of this paper, that the singular locus $S$ of $\F$ is non-empty. Moreover suppose $S$ contains a path-connected component $P$ of type $(A_1)$. It is easy to see that the (topological) closure $\overline{P}$ is also contained in $S$. Hence there exists at least one path-connected component $Q \subset \overline{P} \subset S$ of type $(A_0)$ such that $P$ limits on $Q$. 

Recall from \cite{H2013} that the Legendrian foliation in a tubular neighborhood of $Q$ is well-understood. Namely, a tubular neighborhood $N(Q)$ of $Q \subset W$ is diffeomorphic to a neighborhood of the zero section of a (generally nonlinear) flat $\mathbb{R}^2$-bundle 
\begin{equation*}
	\mathbb{R} ^2 \to E \to Q
\end{equation*}
Denote by $\iota: N(Q) \hookrightarrow E$ the smooth embedding into a neighborhood of the zero section, such that $\iota(Q)$ is the zero section. Let $q \in Q$ be a base point, to be chosen later. Then the flat connection is characterized by the holonomy representation
\begin{equation*}
	\phi: \pi_1(Q,q) \to \textit{Diff}_0^+(\mathbb{R}^2),
\end{equation*}
where $\textit{Diff}_0^+(\mathbb{R}^2)$ is the group of germs of origin-fixing and orientation-preserving diffeomorphisms of $\mathbb{R}^2$ in a neighborhood of the origin. More importantly, it follows from \cite{H2013} that $\iota_\ast(\lambda)$ is a parallel constant vertical one-form on $E$ such that $\iota_\ast(d\lambda)$ restricts to an area form on each fiber $\mathbb{R}^2$. Hereafter we will identify $\lambda$ with $\iota_\ast \lambda$ and will not write $\iota$ anymore. Therefore we have in fact the holonomy respresentation
\begin{equation*}
	\phi: \pi_1(Q,q) \to \mathcal{A}_0^+(\mathbb{R}^2,d\lambda), 
\end{equation*}	
with the image landed in $\mathcal{A}_0^+(\mathbb{R}^2,d\lambda)$, the group of germs of origin-fixing area-preserving diffeomorphisms of $\mathbb{R}^2$ with respect to the area form $d\lambda$.

Now by assumption there is a singular locus $P$ of type $(A_1)$ which limits on $Q$. This implies that there exists a point $q \in Q$, an oriented loop $\gamma \subset Q$ based at $q$, and a sequence of pairwise distinct points $p_k \in E_q \cap P, k=1,2,\dots$ such that the following holds.
\begin{itemize}
	\item $\lambda(p_k)=0$ for any $k$,
	\item $p_{k+1}=\phi_\gamma(p_k)$ for any $k$,
	\item $\lim_{k \to +\infty} p_k=0$,
\end{itemize}
where $\phi_\gamma:=\phi(\gamma) \in \textit{Isom}_0^+(\mathbb{R}^2,d\lambda)$. We will call such $\gamma$ an {\em attracting loop} in $W$.

Now the strategy to prove \autoref{thm:main} is rather simple, namely, we will show that such an attracting loop does not exist. In order to do that, let us briefly recall some useful facts about fixed point in smooth dynamical systems in the next section.

\section{Hyperbolic fixed point theory}

Consider a diffeomorphism $\phi: \mathbb{R}^n \to \mathbb{R}^n$ such that $\phi(0)=0$, we are interested in the local dynamics in a neighborhood of $0$, i.e., the orbit of any point in a sufficiently small neighborhood of the origin under the iterations of $\phi$. The answer to this kind of question belongs to the realm of smooth dynamical systems. One of the most important theorems in this field, and also the one we will make use of, is the celebrated {\em stable manifold theorem}, which we will briefly discuss in this section. In fact, for our purposes, it suffices to consider the case when $n=2$, but it turns out that the particular value of $n$ does not play any role in the following discussions.

We start by setting up some basic terminologies in smooth dynamical systems which will be necessary to state the theorem. Let $\phi: \mathbb{R}^n \to \mathbb{R}^n$ be a diffeomorphism such that $\phi(0)=0$. Suppose $0$ is a {\em hyperbolic fixed point} of $\phi$, i.e., no eigenvalues of $d\phi(0)$ lie on the unit circle. Then there exists a splitting $T_0 \mathbb{R}^n=E^s \oplus E^u$ of vector spaces and a constant $\lambda \in (0,1)$ such that 
\begin{equation*}
	|d\phi(0)(\nu)|<\lambda|\nu| \quad\text{if } \nu \in E^s, \quad\quad |d\phi(0)(\nu)|>\frac{1}{\lambda}|\nu| \quad\text{if } \nu \in E^u,
\end{equation*}
where the norm is taken with respect to some metric on $\mathbb{R}^n$, and it is easy to show that the above definition is independent of the choice of metric. Define the stable and unstable subsets by
\begin{equation*}
	W^s(0)=\{x \in \mathbb{R}^n ~|~ \lim_{k \to +\infty} \phi^k(x)=0\}, \quad\quad W^u(0)=\{x \in \mathbb{R}^n ~|~ \lim_{k \to -\infty} \phi^k(x)=0\}.
\end{equation*}

Now we are in the position to state the following theorem.

\begin{theorem}[Stable manifold theorem] \label{thm:SMT}
Using the notations from above, both $W^s(0)$ and $W^u(0)$ are smooth immersed submanifolds of $\mathbb{R}^n$. Moreover $T_0 W^i(0)=E^i$ where $i \in \{s,u\}$.
\end{theorem}

The possible self-intersections of the (un)stable manifolds are responsible for chaotic dynamical systems such as Smale's horseshoe map, but these phenomena will not bother us since we are only interested in a sufficiently small neighborhood of the origin, where both $W^s(0)$ and $W^u(0)$ are embedded. The interested readers are referred to \cite{Pe1930,Sm1963} and the references therein for the proof of \autoref{thm:SMT} and more thorough treatment.

\section{Proof of \autoref{thm:main}}

The goal of this section is to give a complete proof of \autoref{thm:main}. Using the notations from the end of \autoref{sec:sing}, we assume by contradiction that there exists a singular locus $P$ of type $(A_1)$, which limits on a singular locus $Q$ of type $(A_0)$. Moreover let $\gamma \subset Q$ be a loop based at $q \in Q$ such that $\gamma$ is attracting. Consider the holonomy map $\phi_\gamma: \mathbb{R}^2 \to \mathbb{R}^2$ which fixes the origin. For the rest of the proof we will fix the choice of $\gamma$ and simply write $\phi=\phi_\gamma$. Up to a change of coordinates, we may assume that all the points $p_k \in E_q \cap P$ are contained in the positive $x$-axis such that $p_k \downarrow 0$ as $k \to \infty$. Using the coordinates on $\mathbb{R}^2$, let us write $p_k=(x_k,0)$, and $\phi(x,y)=(\phi^1,\phi^2)$. Since $\phi(0)=0$ and $\phi$ is area-preserving, we have
\begin{equation*}
	\det(d\phi(0))=
	\begin{vmatrix}
		\p_x \phi^1 & \p_y \phi^1 \\
		\p_x \phi^2 & \p_y \phi^2
	\end{vmatrix}_{(x,y)=(0,0)}
	=1. 
\end{equation*}
Moreover it is easy to check, using the fact $\phi(p_k)=p_{k+1}$, that $\p_x \phi^1(0) \leq 1$ and $\p_x \phi^2(0)=0$. We will consider the following two cases separately, namely, either
\begin{equation*}
	d\phi(0)=
	\begin{pmatrix}
		1 & a \\
		0 & 1
	\end{pmatrix}
\end{equation*}
for some $a \in \mathbb{R}$ or $0$ is a hyperbolic fixed point of $\phi$.

Let us first consider the non-hyperbolic case. Let $\text{Id}: \mathbb{R}^2 \to \mathbb{R}^2$ be the identity map. Define $\phi_t=t\text{Id}+(1-t)\phi$ for $0 \leq t \leq 1$. Then $\phi_t$ is also a germ of origin-fixing diffeomorphism of $\mathbb{R}^2$ for any $t$. It is clear that $\phi_t^\ast(\lambda)=\lambda$. This, in particular, implies that $\gamma(\phi_t(p_k))=0$ for any $k$ and $t$. Therefore we obtain a continuous curve
\begin{equation*}
	\tau=\bigcup_{k \gg 0, 0 \leq t \leq 1} \phi_t(p_k) \cup \{0\}
\end{equation*}
containing the origin, such that $\gamma$ vanishes along $\tau$. But this contradicts our assumption that $P$ is a singular locus of type $(A_1)$ because it is contained in a singular locus of type $(A_2).$

Now we turn to the case when $0$ is a hyperbolic fixed point of $\phi$. According to \autoref{thm:SMT}, there exists two smooth curves $W^s(0)$ and $W^u(0)$ passing through the origin. It is clear that $p_k \in W^s(0)$ for all $k$. In the following we will show that $\lambda$ actually vanishes along a (sub-)curve in $W^s(0)$, which will contradict our assumption that $P$ is a singular locus of type $(A_1)$ as before.

Again, up to a change of coordinates, we can assume $W^s(0)$ is contained in the $x$-axis, which covers at least a small neighborhood of the origin. Recall that we can arrange $p_k$, for all $k$, to be contained in the positive $x$-axis and $p_k \downarrow 0$ as $k \to \infty$. In local coordinates, let us write $\lambda=f(x,y)dx+g(x,y)dy$ and $\phi=(\phi^1,\phi^2)$ as before. The equation $\phi^\ast(\lambda)=\lambda$ is equivalent to the following two equations.
\begin{align}
	f &= (f \circ \phi) \p_x \phi^1+(g \circ \phi) \p_x \phi^2 \label{eqn:iterate1} \\
	g &= (f \circ \phi) \p_y \phi^1+(g \circ \phi) \p_y \phi^2 \label{eqn:iterate2}
\end{align}
By restricting to the $x$-axis, we have $\phi^2(x,0) \equiv 0$ by construction. Moreover by the hyperbolicity, we have for sufficiently small $x>0$, $|\p_x \phi^1(x,0)|<1$. Now for sufficiently small $x>0$, (\ref{eqn:iterate1}) implies
\begin{equation*}
	f(x,0)=f(\phi(x,0)) \p_x \phi^1(x,0)=f(\phi^1(x,0),0) \p_x \phi^1(x,0),
\end{equation*}
which can be iterated $n$ times to give
\begin{equation*}
	f(x,0)=f(\phi^{(n)}(x,0)) \prod_{i=1}^n \p_x \phi^1 (\phi^{(i-1)} (x,0)),
\end{equation*}
Letting $n \to \infty$, we see that $f(x,0) \equiv 0$ for sufficiently small $x>0$ since $\lim_{n \to \infty} f(\phi^{(n)} (x,0))=f(0,0)=0$ and $|\p_x \phi_1|$ is bounded from above by 1.

Now we turn to (\ref{eqn:iterate2}), which has the following form when restricted to the $x$-axis
\begin{equation} \label{eqn:iterate3}
	g(x,0)=g(\phi(x,0)) \p_y \phi^2(x,0)
\end{equation}
Recall we write $p_k=(x_k,0)$. Since by construction we have $\phi([x_{k+1},x_k])=[x_{k+2},x_{k+1}]$ being a diffeomorphism for all $k$, the value of $g$ on the positive $x$-axis is completely determined by its value on the interval $[x_2,x_1]$. Here, of course, we assume $x_1>0$ is sufficiently small so that $f$ vanishes on $[0,x_1]$.

Since we will always work on the $x$-axis, so to avoid too many zeros in the upcoming calculations, let us abbreviate notations as follows. We will write
\begin{equation*}
	g(x):=g(x,0), \quad\quad \text{and} \quad \phi^i(x):=\phi^i(x,0) \text{ for } i=1,2.
\end{equation*}
By differentiating both sides of (\ref{eqn:iterate3}) with respect to $x$, we get
\begin{equation} \label{eqn:dotG}
	\p_x g(x)=\p_x g(\phi^1(x)) \p_x \phi^1(x) \p_y \phi^2(x)+g(\phi^1(x)) \p^2_{xy} \phi^2(x).
\end{equation}
If $x>0$ is sufficiently small, so is the second summand in (\ref{eqn:dotG}) compare with the first summand. Suppose there is some $x_0 \in (0,x_1)$ such that $g(x_0) \neq 0$. Then we have the following estimate that $x=x_0$,
\begin{equation} \label{eqn:estimate1}
	|\p_x g(\phi^1({x_0}))|=\left|\frac{\p_x g({x_0})-g(\phi^1({x_0})) \p^2_{xy}\phi^2({x_0})}{\p_x \phi^1({x_0}) \p_y \phi^2({x_0})} \right| \geq (1-\epsilon) \left| \frac{\p_x g({x_0})}{\p_x \phi^1({x_0}) \p_y \phi^2({x_0})} \right|,
\end{equation}
for an arbitrarily small $\epsilon>0$ depending on $x_0$. Now one can iterate (\ref{eqn:estimate1}) to get the following
\begin{align} \label{eqn:estimate2}
	|\p^2_{xx} g(0)| &= \lim_{n \to \infty} \left|\frac{\p_x g(\phi^{1,(n)} (x_0))}{\phi^{1,(n)} (x_0)} \right| \\
					 &\geq \lim_{n \to \infty} (1-\epsilon)^n \left| \frac{\p_x g(x_0)}{\phi^{1,(n)}(x_0) \prod_{i=0}^{n} \p_x \phi^1(\phi^{1,(i)}(x_0) \p_y\phi^2(\phi^{1,(i)} (x_0))} \right| \nonumber\\
					 &\geq \lim_{n \to \infty} (1-\epsilon)^n \left| \frac{\p_x g(x_0)}{x_0 \prod_{i=0}^{n} \big( \lambda\p_x \phi^1(\phi^{1,(i)}(x_0)) \p_y\phi^2(\phi^{1,(i)} (x_0)) \big)} \right| \nonumber\\
					 &= +\infty. \nonumber
\end{align}
Here for the first equality we used the fact that $\p_x g(0)=0$ because $g$ vanishes on a sequence of points $p_k=(x_k,0) \to (0,0)$ as $k \to \infty$. The second inequality is just an iteration of (\ref{eqn:estimate1}), and the third inequality uses our hyperbolicity assumption that $\phi^1(x_0) \leq \lambda x_0$. The last equality uses that fact that $\p_x \phi^1(0) \p_y \phi^2(0) = \det(d\phi(0))=1$ and $\lambda<1$. Of course one should also choose $x_0$ so small such that the factor $(1-\epsilon)^n$ does not effect the overall divergence.

Clearly (\ref{eqn:estimate2}) violates the assumption that $g$ is smooth at 0. Hence $g$ must also vanish on a small interval $[0,\delta]$ in the (positive) $x$-axis, so does $\lambda$. So we see that $P$ is, in fact, contained in a singular locus of type $(A_2)$, which contradicts our original assumption. This finishes the proof of \autoref{thm:main}.

\bibliography{mybib}
\bibliographystyle{amsalpha}

\end{document}